\documentclass[12pt,a4paper]{article}
\usepackage[tbtags]{amsmath}    
\usepackage{amsfonts,amssymb,amsthm,euscript,makeidx,pifont,boxedminipage,multicol,fullpage}
\usepackage{color}

\def\be{\begin{eqnarray}}
\def\ee{\end{eqnarray}}

\def\b*{\begin{eqnarray*}}
\def\e*{\end{eqnarray*}}

\newtheorem{Theorem}{Theorem}[section]
\newtheorem{Lemma}[Theorem]{Lemma}
\newtheorem{Proposition}[Theorem]{Proposition}

\newtheorem{Definition}[Theorem]{Definition}
\newtheorem{Remark}[Theorem]{Remark}
\newtheorem{Example}[Theorem]{Example}

\newtheorem{Assumption}[Theorem]{Assumption}

\newcommand{\rmi}{{\rm (i)$\>\>$}}
\newcommand{\rmii}{{\rm (ii)$\>\>$}}
\newcommand{\rmiii}{{\rm (iii)$\>\>$}}
\newcommand{\rmiv}{{\rm (iv)$\>\>$}}

\def \E{\mathbb{E}}
\def \F{\mathbb{F}}
\def \G{\mathbb{G}}

\def \P{\mathbb{P}}
\def \Q{\mathbb{Q}}
\def \R{\mathbb{R}}

\def\Ac{{\cal A}}

\def\Cc{{\cal C}}
\def\Dc{{\cal D}}

\def\Fc{{\cal F}}
\def\Gc{{\cal G}}

\def\Pc{{\cal P}}
\def\Qc{{\cal Q}}

\def\Sc{{\cal S}}
\def\Tc{{\cal T}}

\def \Om{\Omega}
\def \om{\omega}
\def \Omb{\overline{\Om}}
\def \omb{\bar{\om}}

\def \eps{\varepsilon}

\def \0{\mathbf{0}}
\def \1{\mathbf{1}}
\def \x{\times}
\def \bmu{\boldsymbol{\mu}}

\def \bnu{\boldsymbol{\nu}}

\def \bk{\boldsymbol{K}}
\def \bkk{\boldsymbol{k}}
\def \bc{\boldsymbol{C}}
\def \ba{\boldsymbol{\alpha}}

\def \bl{\boldsymbol{\lambda}}

\def\oc{\overline \Cc}

\def \Fcb{\overline{{\cal F}}}
\def \Fbb{\overline{\F}}

\def \Pcb{\overline{\Pc}}

\def \Pb{\overline{\P}}

\def \Pbb{\mathbf{P}}

\def \db{\overline{d}}

\def \Th{\Theta}

\def \pp {\mathsf{P}}
\def \dd {\mathsf{D}}

\def \no{\noindent}
\def \Gab{\overline{\Gamma}}
\def \SG{\text{SG}}

\title{A stability result on optimal Skorokhod embedding 
	\footnote{Gaoyue Guo is grateful for the helpful discussions he has had with Henry-Labord\`ere, R\"uschendorf, Tan and Touzi, and especially the remarks given by Touzi.}}

\author{
	Gaoyue Guo\thanks{Gaoyue Guo thankfully acknowledges the financial support of the ERC 321111 Rofirm and the ANR Isotace. CMAP, Ecole Polytechnique, France. guo@cmap.polytechnique.fr}
}

\date{\today}

\begin{document}
\bibliographystyle{plain}

\maketitle

\abstract{This is a continuation of Guo, Tan \& Touzi \cite{GTT}. Motivated by the model-independent pricing of derivatives calibrated to the real market, we consider an optimization problem similar to the optimal \textit{Skorokhod embedding problem}, where the embedded Brownian motion needs only to reproduce a finite number of prices of Vanilla options. As same as in \cite{GTT}, we derive in this paper the corresponding dualities and the geometric characterization of optimizers. Then we show a stability result, \textit{i.e.} when more and more Vanilla options are given, the optimization problem converges to an optimal Skorokhod embedding problem, which constitutes the basis of the numerical computation in practice. In addition, by means of different metrics on the space of probability measures, a convergence rate analysis is provided under suitable conditions.}

\vspace{2mm}







\section{Introduction}\label{sec:int}

The Skorokhod embedding problem (SEP) consists in representing a centered probability on the real line as the distribution of a Brownian motion stopped at a chosen stopping time, see \textit{e.g.} the survey paper \cite{Obloj} of Ob{\l}\'oj for a comprehensive account of the field. Motivated by the study of the model-independent pricing of derivatives consistent with the market prices of Vanilla options, the associated optimization problem over all embedding stopping times has received the substantial attention from the mathematical finance community. According to the no-arbitrage framework, the underlying asset is required to be a martingale, and additionally, the market calibration allows to recover the marginal laws of the underlying at certain maturities, see \textit{e.g.} Breeden \& Litzenberger \cite{BL}. Therefore, based on the fact that every continuous martingale is a time-changed Brownian motion by Dambis-Dubins-Schwarz theorem, Hobson studied the robust hedging of lookback options in his seminal paper \cite{Hobson} by means of the SEP. The main idea of his pioneering work is to exploit some solution to the SEP satisfying some optimality criteria, which yields the robust hedging strategy and allows to solve together the model-independent pricing and robust hedging problems. Here after, various extensions were achieved in the literature,  such as Cox \& Hobson \cite{CH}, Hobson \& Klimmek \cite{HK}, Cox, Hobson \& Ob{\l}\'oj \cite{CHO1}, Cox \& Ob{\l}\'oj \cite{CO1},  
Davis, Ob{\l}\'oj \& Raval\cite{DOR} and Ob{\l}\'oj \& Spoida \cite{OS}, \textit{etc}. A thorough literature is provided in Hobson's survey paper \cite{Hobson2011}.

This heuristic idea is generalized by Beiglb\"ock, Cox \& Huesmann \cite{BCH} and Guo, Tan \& Touzi \cite{GTT, GTT3}, where the optimal SEP is defined by a unifying formulation that recovers the previous known results. Namely, their main results are twofold. First, an expected duality is established, \textit{i.e.} identity between the optimal SEP and the corresponding robust hedging problem. Second, they derive respectively a characterization of the optimizers by a geometric pathwise property, due to which, all the known optimal embeddings can be interpreted through one unifying principle.

From the financial viewpoint, Vanilla options are assumed to be ideally liquid in the above literature, \textit{i.e.} the prices of call options are known for all strikes at some maturity, or equivalently, the underlying has the unique distribution determined by the market prices at this maturity. However, only a finite number of call options are traded in practice, which leads to another optimization problem that we shall consider in this paper.

To illustrate the idea, let us start by a simple case where all call options are of the same maturity. Let $\bk=(K_i)_{1\le i\le n}$ be a collection of strikes and $\bc=(C_i)_{1\le i\le n}$ be the set of corresponding prices of call options. Then, roughly speaking, the embedding problem aims to find a stopping time $\tau$ on some Brownian motion $B=(B_t)_{t\ge 0}$ such that the stopped process $B_{\tau \wedge \cdot}:=\big(B_{\tau\wedge t}\big)_{t\ge 0}$ is uniformly integrable (UI) and 
\b*
\E\left[(B_{\tau}-K_i)^+\right]~=~C_i, \mbox{ for all } i=1,\cdots, n.
\e*
The stopping time $\tau$ is called a $(\bk,\bc)-$embedding. In particular, if the price vector $\bc$ is given by some centered probability $\mu$ on $\R$, \textit{i.e.}
\b*
C_i~=~\int (x-K_i)^+d\mu(x), \mbox{ for all } i=1,\cdots, n,
\e*
then every $\mu-$embedding is clearly a $(\bk,\bc)-$embedding. Following the spirit of the optimal SEP, for any given $\bk$ and $\bc$, we may consider an optimization problem among all $(\bk,\bc)-$embeddings, that may write formally as
\b*
\sup_{\tau:~ (\bk,\bc)-\mbox{embedding }}\E\left[\Phi(B,\tau)\right],
\e* 
where $\Phi$ denotes a measurable reward function to be specified later.

Let the no-arbitrage condition hold, namely, there exists a centered distribution $\mu$ that determines the price of each call option. Since the distribution $\mu$ is generally unknown to practitioners, we ask: with the increasing information $\bk^n=(K^n_i)_{1\le i\le n}$ and  $\bc^n=(C^n_i=\int (x-K^n_i)^+d\mu)_{1\le i\le n}$ given by the market, could we recover asymptotically $\mu$? More importantly, does the optimization problem above converge to the optimal SEP with target distribution $\mu$: 
\b*
\sup_{\tau:~ (\bk^n,\bc^n)-\mbox{embedding }}\E\left[\Phi(B,\tau)\right]~ \longrightarrow~ \sup_{\tau:~ \mu-\mbox{embedding }}\E\left[\Phi(B,\tau)\right] \mbox{ as } n~\longrightarrow~\infty.
\e* 
If this convergence holds, we will see later that the computation of $\pp(\mu)$ reduces to a finite-dimensional optimization problem, and moreover, it is of interest to estimate the convergence rate. 

The paper is organized as follows. In Section \ref{sec:oep}, we first recall the formulation of the optimal SEP and then provide the formulation of our optimal embedding problem as well as its dual problems. We establish the two dualities and show further the expected convergence when given more and more market information. Next, the monotonicity principle is derived in Section \ref{sec:mp}. In Section \ref{sec:acr}, by means of different metrics on the space of probability measures, we provide an estimation of the convergence rate. We finally give the related numerical computation in Section \ref{sec:nc}. 

	\vspace{2mm}

	\noindent {\bf Notations.}
	
		\vspace{1mm}
	
	\no \rmi Let $\Om$ be the space of all continuous functions $\om=(\om_t)_{t\ge 0}$ on $\R_+$ such that $\om_0 = 0$, 
	$B=(B_t)_{t\ge 0}$ the canonical process, \textit{i.e.} $B_t(\om):=\om_t$ and $\F= (\Fc_t)_{t \ge 0}$ the canonical filtration generated by $B$, \textit{i.e.} $\Fc_t:=\sigma(B_s, s\le t)$. Denote by $\P_0$ the Wiener measure and by $\F^a= (\Fc^a_t)_{t \ge 0}$ the augmented filtration under $\P_0$.
	
	\vspace{1mm}
	
	\noindent \rmii For every integer $k\ge 1$, define 
	\b*
	\R^k_{<}~ \left(\mbox{resp. } \R^k_>\right)&:=&\left\{(x_1,\cdots, x_k)\in\R^k: x_1> \cdots > x_k~ \left(\mbox{resp. } x_1>\cdots>x_k\right)\right\}
	\e*
	and the corresponding closure $\R^m_{\le}$ (resp. $\R^m_{\ge}$).
		
		\vspace{1mm}

	\noindent \rmiii Let $m \ge 1$ be a fixed integer and set $\Th^m:=\R^m_{\le}~\cap~\R^m_+$. Define the product  canonical space by $\Omb := \Om \x \Th^m$ with its elements denoted by $\omb:=\left(\om,\theta:=(\theta_1,\cdots, \theta_m)\right)$. Denote further by $\left(B, T:=(T_1,\cdots,T_m)\right)$ the canonical element on $\Omb$, \textit{i.e.} $B(\omb):=\om$ and $T(\omb):=\theta$ for all $\omb=(\om,\theta)\in\Omb$. The corresponding canonical filtration is denoted by $\Fbb= (\Fcb_t)_{t \ge 0}$,
	where
	\b*
	\Fcb_t&:=&\sigma\left(B_s, \mbox{ for } s\le t\right)~\vee~\sigma\left(\{T_i\le s\}, \mbox{ for } s\le t \mbox{ and } i=1,\cdots, m\right).
	\e*
	In particular, all random variables $T_1, \cdots, T_m$ are $\Fbb-$stopping times.

	\vspace{1mm}

	\noindent \rmiv We endow
	$\Om$ with the uniform convergence topology,
	and $\Omb$ with the product topology,
	then $\Om$ and $\Omb$ are both Polish spaces (metrizable, separable and complete).


\section{An optimal embedding problem}\label{sec:oep}

In this section, we first give the problem formulation as well as its dual problems, see Guo, Tan \& Touzi \cite{GTT} for more details. We emphasize that the problem is formulated in a weak setting, \textit{i.e.} the stopping times are identified by means of probability measures on the enlarged space $\Omb$.

Let $\Pcb(\Omb)$ be the space of all (Borel)  probability measures on $\Omb$,
	and define
	\be \label{def:st}
		\Pcb
		~:=~
		\Big\{
			\Pb \in \Pcb(\Omb):
			B ~\mbox{is an}~ \Fbb-\mbox{Brownian motion and }
			B_{T_m \wedge \cdot}
			~\mbox{is UI under } \Pb\Big\}.
	\ee 
	\begin{Remark}
	\rmi For any $\Pb\in\Pcb$, $T_1,\cdots, T_m$ can be considered as ordered randomized stopping times under $\Pb$, see \textit{e.g.} Beiglb\"ock, Cox \& Huesmann \cite{BCH}.
	
	\vspace{1mm}

\no \rmii Throughout the paper, $m$ is a fixed integer which stands for the number of maturities. In the following, we consider different embedding problems according to the constraints on $B_{T_1},\cdots, B_{T_m}$.
	\end{Remark}


\subsection{Optimal embedding problems}\label{ssec:oeps}

\paragraph{Optimal Skorokhod embedding problem}

We begin by recalling the optimal SEP. Let $\bmu := (\mu_1, \cdots, \mu_m)$ be a vector of probability distributions on $\R$ and denote, for any integrable function $\phi: \R\to\R$,
\b*
	\mu_i(\phi)~:=~\int_{\R} \phi(x)d\mu_i(x), \mbox{ for all } i=1,\cdots, m.
\e*
The vector $\bmu$ is called a peacock if 
\b*
&&\mu_i \mbox{ has a finite first moment, \textit{i.e.} } \mu_i(|x|)<+\infty, \mbox{ for all } i=1,\cdots, m;  \\
&&\mbox{the map } i\mapsto \mu_i(\phi) \mbox{ is non-decreasing for any convex function } \phi.
\e*
A peacock $\bmu$ is called centered if $\mu_i(x)=0$ for all $i=1, \cdots, m$. Denote by $\Pbb^{\preceq}$ the collection of all centered peacocks. For every vector $\bmu= (\mu_1, \cdots, \mu_m)$, define the set of $\bmu-$embeddings 
\be\label{def:se}
		\Pcb (\bmu)
		&:=&
		\Big\{
			\Pb \in \Pcb:
			B_{T_i}~ \stackrel{\Pb}{\sim}~\mu_i, \mbox{ for all } i = 1, \cdots, m
		\Big\}.
	\ee
As a consequence of Kellerer's theorem, $\Pcb(\bmu)$ is nonempty if and only if $\bmu\in\Pbb^{\preceq}$. In the following,  let us fix a non-anticipative function $\Phi : \Omb \to \R$, \textit{i.e.} $\Phi$ is measurable and $\Phi (\om, \theta)=\Phi \big(\om_{\theta_m \wedge \cdot}, \theta\big)$ holds for all $(\om,\theta)\in \Omb$. We define the optimal SEP for every $\bmu\in\Pbb^{\preceq}$ by
\be \label{def:osep}
		\pp(\bmu)
		&:=&
		\sup_{\Pb \in \Pcb(\bmu)}
		~\E^{\Pb} 
		\left[\Phi(B,T)\right].
	\ee 
It follows by Guo, Tan \& Touzi \cite{GTT} that, under some suitable metric on $\Pbb^{\preceq}$, the map $\bmu\mapsto \pp(\bmu)$ is concave and upper-semicontinuous, which leads to the required dualities. As for the optimal embedding problem in the following, we will proceed with an analogous analysis to derive the dualities.

\paragraph{Optimal embedding problem}

Next let us turn to define the optimization problem described in Section \ref{sec:int}, where the embedded Brownian motion  is only required to be consistent with a finite number of market prices of call options. For the sake of clarity, we assume that the call options of different maturities have the same set of strikes. Throughout the paper, $\bk:=(K_i)_{1\le i\le n}\in\R^n_<$ is reserved for the vector of strikes and $\bc:=(C_{i,j})_{1\le i\le m,1\le j\le n}$ for the price matrix of call options, indexed by maturity and strike. Then a probability $\Pb\in\Pcb$ is called a $(\bk,\bc)-$embedding if 
\be\label{def:em}
\E^{\Pb}\left[(B_{T_i}-K_j)^+\right]~=~C_{i,j}, \mbox{ for all } i=1,\cdots, m \mbox{ and } j=1,\cdots, n.
\ee
Denote by $\Pcb(\bk,\bc)$ the collection of all $(\bk,\bc)-$embeddings and by $\Ac(\bk)\subset\R^{mn}_+$ the set of all matrices $\bc$ such that $\Pcb(\bk,\bc)$ is nonempty. Similarly, for each matrix $\bc\in\Ac(\bk)$, we may define the optimal embedding problem by
\be\label{def:oep}
\pp(\bk,\bc)
		&:=&
		\sup_{\Pb \in \Pcb(\bk,\bc)}
		~\E^{\Pb} 
		\left[\Phi(B,T)\right].
\ee
It follows by definition that $\Ac(\bk)$ is convex, and we pursue to derive the corresponding dualities by a similar argument. However, the set $\Pcb(\bk,\bc)$ is generally not compact, see Example \ref{exa:closure} below. Notice that the map $\bc\mapsto \pp(\bk,\bc)$ is also concave, then a classical result in convex analysis is applied to obtain the required result.
		
\begin{Example}\label{exa:closure}
Take $m=n=1$, $K=0$ and $C=2$. Let $(\mu^k)_{k\ge 2}$ be a sequence of probability distributions defined by
\b*
\mu^k(dx)&:=&\frac{1}{k}\delta_{\{-k\}}(dx)+\left(1-\frac{2}{k}\right)\delta_{\{-\frac{k}{k-2}\}}(dx)+\frac{1}{k}\delta_{\{2k\}}(dx).
\e*
It follows by a straightforward computation that
\b*
\mu^k(x)~=~0 &\mbox{and}& \mu^k(x^+)~=~2.
\e*
Moreover, $\mu^k$ converges weakly to the measure $\mu$ which puts the unit mass on $-1$. Take an arbitrary sequence of measures $(\Pb_k)_{k\ge 2}$ with $\Pb_k\in\Pcb(\mu^k)\subset\Pcb(0,2)$, then it admits a convergent subsequence denoted again by $(\Pb_k)_{k\ge 2}$, see \textit{e.g.} Theorem \ref{theo:conv} or Lemma 4.3 in Guo, Tan \& Touzi \cite{GTT}. Moreover, any accumulation point $\Pb$ satisfies $\E^{\Pb}[B_T]=-1$ and $\E^{\Pb}[B_T^+]=0$, which implies further $\Pb\notin\Pcb(0,2)$.
\end{Example}


\subsection{Dual problems and dualities}\label{ssec:dps}

In this section, we introduce the corresponding dual problems. Recall that $\P_0$ is the Wiener measure on $\Om$ under which the canonical process $B$ is a standard Brownian motion, 
	$\F = (\Fc_t)_{t \ge 0}$ is its natural filtration 
	and $\F^a = (\Fc^a_t)_{t \ge 0}$ is the augmented filtration by $\P_0$.
	Denote by $\Tc^a$ the collection of all increasing families of $\F^a-$stopping times $\tau = (\tau_1, \cdots, \tau_m)$ such that the process $B_{\tau_m \wedge \cdot}$ is uniformly integrable. Define also by $\Lambda$ the space of continuous functions $\lambda:\R\to\R$ with linear growth and by $\Lambda^m$ its $m-$product.
For $\bmu=(\mu_1,\cdots,\mu_m)\in \Pbb^{\preceq}$, $\bl=(\lambda_1,\cdots, \lambda_m)\in\Lambda^m$, $x=(x_1,\cdots, x_m)\in \R^m$ and $\big(\om,\theta=(\theta_1,\cdots,\theta_m)\big)\in\Omb$, we denote
	\b*
	\bmu(\bl)
		~:=~
		\sum_{i =1}^m \mu_i(\lambda_i),~~  
\bl(x)~:=~\sum_{i=1}^m \lambda_i(x_i) ~~\mbox{and}~~ \om_{\theta}~:=~(\om_{\theta_1},\cdots,\om_{\theta_m}).
	\e*
	Then the first dual problem of \eqref{def:osep} is given by
	\be \label{def:dp0}
		\dd_0(\bmu)&:=&\inf_{\bl \in \Lambda^m} \left\{
			\sup_{\Pb \in \Pcb} \E^{\Pb}
			\left[
				\Phi (B,T) 
				-\bl(B_T) 
			\right] 
			+ \bmu(\bl)
		\right\} \\
		&=&
		\inf_{\bl \in \Lambda^m} \left\{
			\sup_{\tau \in \Tc^a} \E^{\P_0}
			\left[
				\Phi (B,\tau) 
				-\bl(B_{\tau}) 
			\right] 
			+ \bmu(\bl)
		\right\}.
	\ee
As for the second dual problem, we return to the enlarged space $\Omb$. An $\Fbb-$adapted continuous process $S= (S_t)_{t \ge 0}$ is called a $\Pcb-$strong supermartingale if
\b*
		\E^{\Pb} \left[ S_{\tau_2}\big |\Fcb_{\tau_1} \right]
		~\le~
		S_{\tau_1},
		~\Pb-\mbox{a.s.} 
	\e*	
	for all $\Fbb-$stopping times $\tau_1 \le \tau_2$ and all $\Pb\in\Pcb$. Denote by $\Sc$ the set of all $\Pcb-$strong supermartingales starting at zero and put
\b*
		\Dc &:=&
		\left\{ 
			(\bl, S) \in \Lambda^m \x \Sc
			:
			\bl(\om_{\theta}) +S_{\theta_m}(\om)
			~\ge~ \Phi (\om,\theta),~ \mbox{for all } (\om,\theta)\in\Omb
		\right\}.
	\e*
Then the second dual problem is given by
	\be\label{def:dp}
		\dd(\bmu)
		&:=&
		\inf_{(\bl, S) \in \Dc}\bmu(\bl).
	\ee
	\begin{Remark}
\rmi By penalizing the marginal constraints, we obtain the first dual problem $\dd_0(\bmu)$, where a multi-period optimal stopping problem appears for every fixed $\bl \in \Lambda^m$. 

	\vspace{1mm}

\no \rmii The second dual problem $\dd(\bmu)$ is slightly different to that in \cite{GTT}, where we wrote $S$ as an stochastic integral. As our main concern is to study the model-independent pricing, we show later that the formulation $\dd(\bmu)$ is enough to deduce the required results.  
\end{Remark}

\no Following the idea above, we may define similarly the dual problems of \eqref{def:oep}. For $\bc=(C_{i,j})_{1\le i\le m,1\le j\le n}$, $\ba=(\alpha_{i,j})_{1\le i\le m, 1\le j\le n}\in\R^{mn}$ and $x=(x_1,\cdots, x_n)\in\R^n$, set
\b*
 \ba\cdot \bc~:=~\sum_{i=1, j=1}^{m,n}\alpha_{i,j}\cdot C_{i,j},~~ \ba\cdot x~:=~\sum_{i=1, j=1}^{m,n}\alpha_{i,j}\cdot x_j~~\mbox{and}~~x^+~:=~(x_1^+,\cdots, x_n^+). 
\e*
Hence, the dual problems of \eqref{def:oep} are defined respectively by
\be \label{def:adp0}
		\dd_0(\bk,\bc)&:=&\inf_{\ba \in \R^{mn}} \Big\{
			\sup_{\Pb \in \Pcb} \E^{\Pb}
			\big[
				\Phi (B,T) - 
				\ba\cdot (B_{T}-\bk)^+ 
			\big] 
			+ \ba\cdot \bc
		\Big\} \\
		&=&
		\inf_{\ba \in \R^{mn}} \Big\{
			\sup_{\tau \in \Tc^a} \E^{\P_0}
			\big[
				\Phi (B,\tau) 
				-\ba\cdot (B_{\tau}-\bk)^+ 
			\big] 
			+  \ba\cdot \bc
		\Big\}.
	\ee
	and
	\be\label{def:adp}
		\dd(\bk,\bc)
		&:=&
		\inf_{(\ba, S) \in \Dc(\bk)}\ba\cdot \bc.
	\ee
where
\b*
		\Dc(\bk) ~:=~
		\left\{ 
			(\ba, S) \in \R^{mn} \x \Sc
			:
			\ba\cdot (\om_{\theta}-\bk)^+  +S_{\theta_m}(\om)
			~\ge~ \Phi (\om,\theta),~ \mbox{ for all } (\om,\theta)\in\Omb
		\right\}.
	\e*
	
\begin{Assumption} \label{hyp:phi1}
		$\Phi:\Omb\to\R$ is bounded from above and the map $\omb \mapsto \Phi(\omb)$ is upper-semicontinuous.	\end{Assumption}
	
	\begin{Theorem} \label{theo:duality}
		Let Assumption \ref{hyp:phi1} hold. 
		
		\no \rmi There exists a $\Pb^* \in \Pcb(\bmu)$ such that
		\b* 
		\E^{\Pb^*}\left[\Phi(B,T)\right] ~~=~~ \pp(\bmu) ~~=~~ \dd_0 (\bmu)~~=~~\dd(\bmu).
		\e*
		\rmii  Assume further $\bc\in\Cc(\bk)$, then one has 
		\b*
		\pp(\bk,\bc)~~=~~\dd_0(\bk,\bc)~~=~~\dd(\bk,\bc).
		\e*
	\end{Theorem}
	
	\proof The proof of $\dd_0(\bk,\bc)=\dd(\bk,\bc)$ is as same as that of $\dd_0(\bmu)=\dd(\bmu)$ by by Guo, Tan \& Touzi \cite{GTT}, then It remains to prove $\pp(\bk,\bc)=\dd_0(\bk,\bc)$. 
	
	One has by Corollary 4.2 of Davis \& Hobson \cite{DH} that $\Cc(\bk)\subset\Ac(\bk)\subset\oc(\bk)$, where $\oc(\bk)$ denotes the closure of $\Cc(\bk)$ and $\Cc(\bk)\subset\R^{mn}_+$ consists of all matrices $\bc=(C_{i,j})_{1\le i\le m,1\le j\le n}$ such that for all $i=1, \cdots, m$ and $j=1, \cdots, n$
\b*
(C_{i,j})_{1\le i\le m}~\in~\R^m_< &\mbox{and}&  (C_{i,j})_{1\le j\le n} ~\in~ \R^n_>, \\
C_{i,j}~>~(-K_j)^+ &\mbox{and}& 1~>~ \frac{C_{i,j}-C_{i,j+1}}{K_{j+1}-K_j}~>~ \frac{C_{i,j-1}-C_{i,j}}{K_{j}-K_{j-1}}. 
\e*
	
	For the sake of simplicity, we set $\pp(\bc)\equiv \pp(\bk,\bc)$ throughout the proof. It follows by definition that the map $\bc\mapsto \pp(\bc)$ is concave on $\Cc(\bk)$. Notice that $\Cc(\bk)\subset \R^{mn}$ is convex and open, then $\bc\mapsto \pp(\bc)$ is continuous on $\Cc(\bk)$. Hence, we may follow the reasoning in the proof of Theorem 3.10 in Guo, Tan \& Touzi \cite{GTT2} and extend the map $\pp$ from $\Cc(\bk)$ to $\R^{mn}$ by
	\b*
		\widetilde{\pp}(\bc)
		&:=&
		\begin{cases}
			 \pp(\bc), ~~&\mbox{ if } \bc\in\Cc(\bk), \\
			-\infty, ~~& \mbox{ otherwise}.
		\end{cases}
	\e* 
 	The concavity and upper semicontinuity of the map $\bc \mapsto \widetilde{\pp}(\bc)$ follow immediately from the definition. Then, one obtains by Fenchel-Moreau theorem
	\b*
		\widetilde{\pp}(\bc)&=&\widetilde{\pp}^{\ast\ast}(\bc),
	\e*
	where $\widetilde{\pp}^{\ast\ast}$ denotes the biconjugae of $\widetilde{\pp}$. In particular, for $\bc\in\Cc(\bk)$ one has
	\b*
		&&
		\pp(\bc) ~~=~~ \widetilde{\pp}(\bc) ~~=~~ \widetilde{\pp}^{\ast\ast}(\bc)\\
		&=&
		\inf_{\ba\in\R^{mn}}\big\{\ba\cdot \bc-\widetilde{\pp}^{\ast}(\ba)\big\}
		~~=~~
		\inf_{\ba\in\R^{mn}}\Big\{\ba\cdot \bc-\inf_{\bc'\in\R^{mn}}\big\{\ba\cdot \bc'-\widetilde{P}(\bc')\big\}\Big\} \\
		&\ge&
		\inf_{\ba\in\R^{mn}}
		\Big\{\ba\cdot \bc
			+
			\sup_{\bc\in\Cc(\bk)} \Big\{\sup_{\P\in\Pcb(\bk,\bc)} \E^{\P} \big[\Phi (B,T) - 
				\ba\cdot (B_{T}-\bk)^+ \big] \Big\}
		\Big\} \\
		&=&
		\inf_{\ba\in\R^{mn}}
		\Big\{\ba\cdot \bc
			+
			\sup_{\P\in\cup_{\bc\in\Cc(\bk)}\Pcb(\bk,\bc)}\Big\{ \E^{\P} \big[\Phi (B,T) - 
				\ba\cdot (B_{T}-\bk)^+ \big] \Big\}
		\Big\} \\
		&=&
		\inf_{\ba\in\R^{mn}}
		\Big\{\ba\cdot \bc
			+
			\sup_{\P\in\cup_{\bc\in\oc(\bk)}\Pcb(\bk,\bc)}\Big\{ \E^{\P} \big[\Phi (B,T) - 
				\ba\cdot (B_{T}-\bk)^+ \big] \Big\}
		\Big\},
		\e*
		where the last inequality follows from the fact that $\cup_{\bc\in\oc(\bk)}\Pcb(\bk,\bc)$ is included in the closure of $\cup_{\bc\in\Cc(\bk)}\Pcb(\bk,\bc)$ under the weak convergence. Hence
		\b*
	\pp(\bc) ~~=~~ \widetilde{\pp}(\bc) ~~=~~ \widetilde{\pp}^{\ast\ast}(\bc)~~\ge~~\dd_0(\bc)~~\ge~~\pp(\bc),
		\e*
		which yields the required duality.
	\qed
	
\begin{Remark}\label{rem:existence}
Notice that the set $\Pcb(\bk,\bc)$ is generally not compact with respect to the weak convergence, due to which the existence of optimizers of $\pp(\bk,\bc)$ can not be ensured. 
\end{Remark}


\subsection{Convergence of optimal embedding problems}\label{ssec:sta}

Let us study here the asymptotic behavior of the upper bound with respect to the market information. Assume that the market is consistent with a centered peacock $\bmu=(\mu_1,\cdots,\mu_m)$, then we ask: when more and more call options are traded, does the upper bound converge to $\pp(\bmu)$? Namely, let $\bk^n=(K^n_i)_{1\le i\le n}\in\R^n_<$ be the vector of strikes that are available in the market and $\bc^n=(C^n_{i,j})_{1\le i\le m, 1\le j\le n}$ be the corresponding price matrix given by $C^n_{i,j}=\mu_i\left((x-K^n_j)^+\right)$ for all $i=1,\cdots, m$ and $j=1,\cdots, n$. Define the bound and the mesh of $\bk$ by
\b*
|\bk^n|~:=~(K_1^n)^-\wedge (K^n_n)^+ &\mbox{and}& \Delta\bk^n~:=~\max_{1<i\le n}\left(K^n_i-K^n_{i-1}\right).
\e*
Loosely speaking, to capture asymptotically $\bmu$ by $(\bk^n,\bc^n)$, a necessary condition for the sequence $(\bk^n)_{n\ge 1}$ is the following:
\begin{Assumption} \label{hyp:stab}
The sequence $(\bk^n)_{n\ge 1}$ satisfies 
\b*
\lim_{n\to\infty}|\bk^n|~=~+\infty&\mbox{and}&\lim_{n\to\infty}\Delta \bk^n~=~0.
 \e*
\end{Assumption} 

\no Motivated by financial applications, the sequence $(\bk^n)_{n\ge 1}$ is assumed to be increasing, \textit{i.e.} $\bk^n\subset\bk^{n+1}$ for all $n\ge 1$ if $\bk^n$ and $\bk^{n+1}$ are viewed as sets. It follows by definition that the map $n\mapsto\pp(\bk^n,\bc^n)$ is non-increasing and 
\b*
\pp(\bk^n,\bc^n)~\ge~ \pp(\bmu),~ \mbox{for all } n\ge 1.
\e*

\begin{Theorem}\label{theo:conv}
Let Assumption \ref{hyp:phi1} hold. Then for any increasing sequence $(\bk^n)_{n\ge 1}$ satisfying Assumption \ref{hyp:stab}, one has
\b*
\lim_{n\to\infty}\pp(\bk^n,\bc^n)&=&\pp(\bmu),
\e*
where $\bc^n$ is defined by $\bmu$ as above.
\end{Theorem}

\proof Notice that the sequence $\pp(\bk^n,\bc^n)$ is non-increasing and thus the limit exists. Let $(\Pb_n)_{n\ge 1}$ be a sequence such that $\Pb_n\in\Pcb(\bk^n,\bc^n)$ and
\b*
\lim_{n\to\infty}\pp(\bk^n,\bc^n)&=&\lim_{n\to\infty}\E^{\Pb_n}\left[\Phi(B,T)\right].
\e*
Then repeating the reasoning of Lemma 4.3 in Guo, Tan \& Touzi \cite{GTT}, we deduce that the sequence $(\Pb_n)_{n\ge 1}$ is tight and any accumulation point $\Pb$ of $(\Pb_n)_{n\ge 1}$ belongs to $\Pcb$. Without loss of generality, denote again by $(\Pb_n)_{n\ge 1}$ the convergent subsequence with limit $\Pb$, then it follows by Lemma \ref{lem:stab} that 
\b*
B_{T_i}~\stackrel{\Pb}{\sim}~\mu_i, \mbox{ for all } i=1,\cdots, m,
\e*
which implies that $\Pb\in\Pcb(\bmu)$. The proof is fulfilled by Fatou's lemma:
\b*
\pp(\bmu)~\le~\lim_{n\to\infty}\pp(\bk^n,\bc^n)~=~\lim_{n\to\infty}\E^{\Pb_n}\left[\Phi(B,T)\right]~\le~\E^{\Pb}\left[\Phi(B,T)\right]~\le~\pp(\bmu).
\e*
 \qed

\begin{Remark}
Combining Theorems \ref{theo:conv} and \ref{theo:duality}, the $\bl$ appearing in $\dd_0(\bmu)$ and $\dd(\bmu)$ can be restricted to take values in the set of Lipschitz functions.
\end{Remark}

\begin{Lemma}\label{lem:stab}
Let $\mu$ be a probability measure on $\mathbb R$ and set $c(K):=\mu\left((x-K)^+\right)$ for every $K\in\R$. Let $(\mu^n)_{n\ge 1}$ be a weakly convergent sequence of probability measures such that
\b*
\mu^n\left(\left(x-K^n_i\right)^+\right)~=~\mu\left((x-K^n_i)^+\right), \mbox{ for all } i=1,\cdots, n,
\e*
where $\left(\bk^n=(K^n_i)_{1\le i\le n}\right)_{n\ge 1}$ is an increasing sequence satisfying Assumption \ref{hyp:stab}. Then 
\b*
\lim_{n\to\infty}\mu^n&=&\mu.
\e*
\end{Lemma}

\proof Set $c^n(K):=\mu^n\left((x-K)^+\right)$ (resp. $c(K):=\mu\left((x-K)^+\right)$) for all $K\in\R$. Since $\left|(x-a)^+ - (x-b)^+\right| \le |a-b|$, the function $K\mapsto c^n(K)$  (resp. $K\mapsto c(K)$) is $1-$Lipschitz. Moreover, for every $K\in \cup_{n\ge 1}\bk^n$, one has $c^n(K)=c(K)$ for $n$ large enough, which implies that $c^n$ converges uniformly to $c$ as $\cup_{n\ge 1}\bk^n$ is dense on $\R$.  

Let $f:\R\to\R$ be twice differentiable with compact support, then it follows by Carr-Madan's formula that
\b*
f(x)&=&\int_{\R} f''(K)(x-K)^+dK.
\e*
We obtain in view of Fubini's Theorem 
\b*
\mu^n(f) &=& \int_{\R} f''(K) c^n(K) dK,
\e* 
which implies that $\mu^n(f)\rightarrow \mu(f)$ as $n\to\infty$. As any continuous function with compact support can be uniformly approximated by a twice differentiable function with compact support, one has $\lim_{n\to\infty}\mu^n=\mu$. \qed


\subsection{Additional market information: Power option}

As described in Remark \ref{rem:existence}, the set of embeddings $\Pcb(\bk,\bc)$ is generally not compact. For technical reasons, we consider in the following two subsets $\Pcb_V$ and $\Pcb_V(\bk,\bc)$:
\b*
\Pcb_V~:=~\left\{\Pb\in\Pcb: \E^{\Pb}\left[|B_{T_m}|^p\right]~=~V\right\} &\mbox{and}& \Pcb_V(\bk,\bc)~:=~\Pcb_V~\cap~ \Pcb(\bk,\bc),
\e*
where $p>1$ and $V$ are fixed throughout the paper. The restriction of the embeddings comes from a new information: the Power option of the last maturity is observed in the market. Consequently, this implies the unknown peacock $\bmu$ must satisfy
\b*
\mu_m(|x|^p)&=&V.
\e*
Put similarly
\b*
\pp^V(\bk,\bc)&:=&\sup_{\Pb \in \Pcb_V(\bk,\bc)}
		~\E^{\Pb} 
		\left[\Phi(B,T)\right].
\e*
and
\b*
		\dd_0^V(\bk,\bc)&:=&\inf_{\ba \in \R^{mn}} \Big\{
			\sup_{\Pb \in \Pcb_V} \E^{\Pb}
			\big[
				\Phi (B,T) - 
				\ba\cdot (B_{T}-\bk)^+ 
			\big] 
			+ \ba\cdot \bc
		\Big\} \\
		&=&
		\inf_{\ba \in \R^{mn}} \Big\{
			\sup_{\tau \in \Tc_V^a} \E^{\P_0}
			\big[
				\Phi (B,\tau) 
				-\ba\cdot (B_{\tau}-\bk)^+ 
			\big] 
			+  \ba\cdot \bc
		\Big\},
\e*
where $\Tc_V^a$ denotes the subset of $\Tc^a$ consisting of elements $\tau$ such that $\E^{\P_0}[|B_{\tau_m}|^p]=V$. Then using exactly the same arguments in Guo, Tan \& Touzi \cite{GTT}, one has the following theorem.

\begin{Theorem} \label{theo:duality2}
Let Assumption \ref{hyp:phi1} hold. Then for any $\bc\in\Ac(\bk)$, there exists a $\Pb^* \in \Pcb_V(\bk,\bc)$ such that
		\b*
		\E^{\Pb^*}\left[\Phi(B,T)\right] ~~=~~\pp^V(\bk, \bc) ~~=~~ \dd^V_0 (\bk,\bc).
		\e*
\end{Theorem}


\section{Monotonicity principle: $m=1$}\label{sec:mp}

Similar to the monotonicity principle introduced in Beiglb\"ock, Cox \& Huesmann \cite{BCH}, we give on one-marginal case  another principle that links the optimality of an embedding and the geometry of its support set. 

For every $\omb=(\om,\theta), \omb'=(\om',\theta')\in\Omb$,
	we define the concatenation $\omb\otimes\omb'\in\Omb$ by
	\b*
		\omb\otimes\omb'&:=&(\om\otimes_{\theta}\om', \theta+\theta'),
	\e* 
	where $\om\otimes_{\theta}\om'\in\Om$ is defined by
	\b*
		\left(\om\otimes_{\theta}\om'\right)_t
		~:=~
		\om_t\mathbf{1}_{[0,\theta)}(t)
		+
		\left(\om_{\theta}+\om'_{t-\theta}\right)\mathbf{1}_{[\theta,+\infty)}(t),
		~\mbox{for all}~ t\in\R_+.
	\e*	
Let $\Gab\subseteq\Omb$ be a subset, we define $\Gab^<$ by
	\b*
		\Gab^<
		&:=&
		\left\{\omb=(\om,\theta)\in\Omb
			~:
			\omb = \omb'_{\theta \wedge \cdot}
			~\mbox{for some}~
			\omb' \in \Gab
			~\mbox{with}~
			\theta' > \theta
		\right\}.
	\e*
	\begin{Definition} \label{def:SG}
		A pair $(\omb,\omb')\in \Omb\x\Omb$ is said to be a stop-go pair relative if 
		$\om_{\theta}=\om'_{\theta'}$ and 
		\b*
			\xi(\omb)+\xi(\omb'\otimes \omb'')~>~\xi(\omb\otimes \omb'')+\xi(\omb')
			&\mbox{for all}&
			\omb'' \in\Omb^+,
		\e*
		where $\Omb^{+}:=\big\{\omb=(\om,\theta)\in\Omb: \theta>0\big\}$.
		Denote by $\SG$ the set of all stop-go pairs.
	\end{Definition}
	
\no By exactly the same arguments in Guo, Tan \& Touzi \cite{GTT3}, we have the following theorem.	
	
	\begin{Theorem}\label{Th:MP}
		Suppose that the optimal embedding problem $\pp(\bk,\bc)$ admits an optimizer 
		$\Pb^{\ast}\in\Pcb(\bk,\bc)$, and the duality $\pp(\bk,\bc) = \dd(\bk,\bc)$ holds true.		
		Then there exists a Borel subset $\Gab^{\ast}\subseteq\Omb$ such that
		\be \label{eq:GC}
			\Pb^* \big[ \Gab^* \big] = 1
			~&\mbox{and}&~
			\SG ~\cap~\big(\Gab^{\ast<}\x\Gab^{\ast}\big)
			~=~
			\emptyset.
		\ee
	\end{Theorem}


\section{Analysis of convergence rate}\label{sec:acr}

As we have shown in Section \ref{ssec:sta} the convergence of $\pp(\bk^n,\bc^n)$ to $\pp(\bmu)$ for any increasing sequence $(\bk^n)_{n\ge 1}$ satisfying Assumption \ref{hyp:stab}, we continue to estimate the convergence rate in this section. First, notice by definition that
\b*
\pp(\bmu)~~\le~~\pp^V(\bk^n,\bc^n)~~\le~~\pp(\bk^n,\bc^n),~ \mbox{ for all } n\ge 1,
\e*
which implies by Theorem \ref{theo:conv} that
\b*
\lim_{n\to\infty}\pp^V(\bk^n,\bc^n)&=&\pp(\bmu).
\e*
Throughout this section we focus on the asymptotic behavior $\pp^V(\bk^n,\bc^n)$ by assuming that $V=\mu_m(|x|^p)<+\infty$, where $p>1$ is a fixed number and $q>1$ is defined by $1/p+1/q=1$. 


\subsection{Metrics on probability space}

In preparation, let us introduce two metrics that are used in the following. Denote  respectively by $\rho(\cdot,\cdot)$ the L\'evy-Prokhorov metric and by $W(\cdot,\cdot)$ the Wasserstein metric, \textit{i.e.} for any two probability measures $\mu$ and $\nu$  on $\R$, one has
\b*
\rho(\mu,\nu)&:=&\inf\left\{\eps>0: F_{\mu}(x-\eps)-\eps\le F_{\nu}(x)\le F_{\mu}(x+\eps)+\eps \mbox{ for all } x\in\R\right\}, \\
W(\mu,\nu)&:=&\inf_{\chi\in\Pc(\mu,\nu)}\int_{\R^2}|x-y|d\chi(x,y),
\e* 
where $F_{\mu}$ (resp. $F_{\nu}$) denotes the cumulative distribution function of $\mu$ (resp. $\nu$), and $\Pc(\mu,\nu)$ denotes the collection of probability measures on $\R^2$ with marginal distributions $\mu$ and $\nu$.

\begin{Remark}\label{rem:rate}
It is well known that, see \textit{e.g.} Chapter 1.2 of Rachev \& R\"uschendorf, \cite{RR}
\b*
W(\mu,\nu)&=&\int_{\R}\left|F_{\mu}(x)-F_{\nu}(x)\right|dx \\
&=&\inf\left\{\mu(f)-\nu(f): f: \R\to\R \mbox{ is } 1-\mbox{Lipschitz}\right\}.
\e*
\end{Remark}

\begin{Lemma}\label{lem:rate}
Let $\mu$ and $\nu$ be two probability measures supported on $[-R,R]$ for some fixed $R>0$. Assume that there exists  some $\eps>0$ such that
\b*
\left|\mu\left((x-K)^+\right)-\nu\left((x-K)^+\right)\right|~\le~\eps \mbox{ for all } K\in [-R,R].
\e*
Then
\b*
\rho(\mu,\nu)~\le~\sqrt{2\eps} &\mbox{and}& W(\mu,\nu)~\le~4R\sqrt{\eps}.
\e*
\end{Lemma}

\proof \rmi Take an arbitrary $0<\delta<\rho(\mu,\nu)$, then one has some $K\in [-R,R]$ such that $F_{\mu}(K-\delta)-F_{\nu}(K)>\delta$ or $F_{\nu}(K)-F_{\mu}(K+\delta)>\delta$. Take the first case without loss of generality, which yields
\b*
\int_{K-\delta}^K\left(F_{\mu}(x)-F_{\nu}(x)\right)dx~\ge~\int_{K-\delta}^K\left(F_{\mu}(K-\delta)-F_{\mu}(K)\right)dx~>~\delta^2.  
\e*
In addition,
\b*
&&\int_{K-\delta}^K\left(F_{\mu}(x)-F_{\nu}(x)\right)dx \\
&=&\left|\int_{K-\delta}^{+\infty}\left(F_{\mu}(x)-F_{\nu}(x)\right)dx - \int^{+\infty}_{K}\left(F_{\mu}(x)-F_{\nu}(x)\right)dx\right| \\
&=&\left|\mu\left((x-K+\delta)^+\right)-\nu\left((x-K+\delta)^+\right)-\mu\left((x-K)^+\right)+\nu\left((x-K+\delta)^+\right)\right|.\e*
It follows by assumption that $\delta^2<2\eps$. That is, 
\b*
\rho(\mu,\nu)&\le&\sqrt{2\eps}.
\e*
\rmii It follows by Theorem 1.1.8 in Rachev \& R\"uschendorf \cite{RR}, see also Remark \ref{rem:rate}, that 
\b*
W(\mu,\nu)&=&\int_{\R}\left|F_{\mu}(x)-F_{\nu}(x)\right|dx \\
&=&\int_{-R}^R\left|F_{\mu}(x)-F_{\nu}(x)\right|dx.
\e*
In addition, it follows by definition that $|F_{\mu}(x)-F_{\nu}(x)|\le \sqrt{2}\rho(\mu,\nu)\le 2\sqrt{\eps}$ for all $x\in [-R,R]$. Thus
\b*
W(\mu,\nu)&\le&4R\sqrt{\eps}.
\e*
\qed

\begin{Proposition}\label{pro:rate}
For each $n\ge 1$ and any $\Pb\in\Pcb(\bk^n, \bc^n)$, set $\nu_i=\Pb~\circ~ (B_{T_i})^{-1}$ for all $i=1,\cdots, m$. Then there exists a constant $C>0$ depending only on $V$ such that, for all $i=1,\cdots, m$ one has
\b*
\rho(\mu,\nu)~\le~C\left(\sqrt{\Delta \bk^n}+|\bk^n|^{-p/2q}\right) \mbox{ and } W(\mu,\nu)~\le~C|\bk^n|\left(\sqrt{\Delta \bk^n}+|\bk^n|^{-p/2q}\right).
\e*
\end{Proposition}

\proof Without loss of generality, it suffices to show the inequality for $i=m$. For the sake of simplicity, we write $\mu\equiv \mu_m$ and $\nu\equiv \nu_m$. The main idea is to approximate $\mu$ and $\nu$ by their truncated versions. 

Set $R=|\bk^n|$ and let $\mu'$ (resp. $\nu'$) be the truncated distribution of $\mu$ (resp. $\nu$). Indeed, let $X$ (resp. $Y$) denotes some random variable of law $\mu$ (resp. $\nu$), then $\mu'$ (resp. $\nu'$) be the law of $X':=(-R)\vee(R\wedge X)$ (resp. $Y':=(-R)\vee(R\wedge Y)$). Then one has for all $K\in [-R,R]$
\b*
&&\left|\mu'\left((x-K)^+\right)-\nu'\left((x-K)^+\right)\right|~=~\left|\E\left[(X'-K)^+\right]-\E\left[(Y'-K)^+\right]\right| \\
&\le&\left|\E\left[(X-K)^+\right]-\E\left[(Y-K)^+\right]\right| +\E\left[|X-X'|\right] +\E\left[|Y-Y'|\right] \\
&\le&\left|\mu\left((x-K)^+\right)-\nu\left((x-K)^+\right)\right| +2\E\left[|X|\mathbf{1}_{\{|X|> R\}}\right]+2\E\left[|Y|\mathbf{1}_{\{|Y|> R\}}\right].~~ (\ast)  \\
\e*
Notice that for each $K\in [-R,R]$ there exists $1\le i<n$ such that $K\in [K^n_i,K^n_{i+1}]$, thus
\b*
\mu\left((x-K)^+\right)-\nu\left((x-K)^+\right)&\le&\mu\left((x-K_i^n)^+\right)-\nu\left((x-K_{i+1}^n)^+\right) \\
&=&\mu\left((x-K_i^n)^+\right)-\mu\left((x-K_{i+1}^n)^+\right)~~\le~~\Delta \bk^n. \\
\e*
Hence
\b*
\left|\mu\left((x-K)^+\right)-\nu\left((x-K)^+\right)\right|&\le& \Delta \bk^n.
\e*
In addition,
\b*
\E\left[|X|\mathbf{1}_{\{|X|> R\}}+|Y|\mathbf{1}_{\{|Y|> R\}}\right]&\le&\E\left[\frac{|X|^p}{R^{p-1}}\mathbf{1}_{\{|X|> R\}}+\frac{|Y|^p}{R^{p-1}}\mathbf{1}_{\{|Y|> R\}}\right]~~\le~~\frac{2V}{R^{p-1}},
\e*
which yields by $(\ast)$ that
\b*
\left|\mu'\left((x-K)^+\right)-\nu'\left((x-K)^+\right)\right|&\le&\Delta\bk^n+\frac{4V}{|\bk^n|^{p/q}}.
\e*
It follows by Lemma \ref{lem:rate} that there exists some $C>0$ such that 
\b*
\rho(\mu',\nu')~\le~C\left(\sqrt{\Delta \bk^n}+|\bk^n|^{-p/2q}\right) \mbox{ and } W(\mu',\nu')~\le~C|\bk^n|\left(\sqrt{\Delta \bk^n}+|\bk^n|^{-p/2q}\right).
\e*
It remains to estimate $\rho(\mu,\mu')$ (resp. $\rho(\nu,\nu')$) and $W(\mu,\mu')$ (resp. $W(\nu,\nu')$). It follows by definition
\b*
\rho(\mu,\mu')~ (\mbox{resp. } \rho(\nu,\nu'))&\le&\frac{2V}{|\bk^n|^p}, \\
W(\mu,\mu')~ (\mbox{resp. } W(\nu,\nu'))&\le&\frac{2V}{|\bk^n|^{p/q}},
\e*
which yield the required inequalities by the triangle inequality. \qed

\begin{Remark}
Notice that, in order to ensure $W(\mu,\nu)$ converge to zero, we need $p>3$ and $|\bk^n|\sqrt{\Delta \bk^n}\to 0$ as $n\to\infty$. 
\end{Remark}

\no To estimate the convergence, we need more regularity on $\Phi$. Let us formulate the assumption on $\Phi$. First, we introduce a metric $\db$ on $\Omb$. 
	for all  $\omb=(\om,\theta_1,\cdots, \theta_m)$, $\omb'=(\om',\theta_1',\cdots, \theta_m')\in\Omb$,
	\b*
		\db(\omb,\omb')
		&:=&
		\sum_{i= 1}^m \left(
			\sqrt{|\theta_i-\theta_i'|}
			+
		|\om_{\theta_i\wedge\cdot}-\om'_{\theta_i'\wedge\cdot}|
		\right),
	\e*
	where
	\b*
		|\om_{\theta_i\wedge\cdot}-\om'_{\theta_i'\wedge\cdot}|
&:=&\sup_{t\ge 0}|\om_{\theta_i\wedge t}-\om'_{\theta_i'\wedge t}|.
	\e*
	We end this section by the following assumption on $\Phi$:
	
\begin{Assumption} \label{hyp:phi0}
$\Phi:\Omb\to\R$ is uniformly bounded and $\db-$Lipschitz, where the uniform norm and Lipschitz constant are respectively denoted by $\|\Phi\|$ and $L$.
	\end{Assumption}


\subsection{One-marginal case: $m=1$}

We start by the one-marginal case, where we may construct explicitly an approximation of a given martingale. Through this, we obtain the difference of the upper bounds between different target distributions under the following assumption.

\begin{Assumption} \label{hyp:phi-one}
$\Phi$ is time-invariant, that is, $\Phi(\om,\theta)=\Phi\left(\om_{\varphi},\varphi(\theta)\right)$ holds for any $(\om,\theta)\in\Omb$ and any increasing function $\varphi: \R_+\to\R_+$.
\end{Assumption}

\begin{Proposition}\label{pro:one-rate}
Let Assumptions \ref{hyp:phi0} and \ref{hyp:phi-one} hold. Then there exists a constant $C>0$ such that, for any centered probability distribution $\nu$ on $\R$ satisfying $\nu(|x|^p)\le V$, one has 
\b*
\left|\pp(\mu)-\pp(\nu)\right|&\le &C\rho(\mu,\nu)^{1/2q}.
\e*
\end{Proposition}

\proof Write $\rho:= \rho(\mu,\nu)$ for the sake of simplicity. Take a $\Pb\in\Pcb(\mu)$, then one has by definition $B_T\stackrel{\Pb}{\sim}\mu$. It follows by Theorem 4 on page 358 in Shiryaev \cite{Shiryaev} and Theorem 1 in Skorokhod \cite{Skorokhod} that, there exist a measurable function $f: \R^2\to\R$ and a Gaussian random variable $G$ that is independent of $\Fbb$ such that
\b*
M~:=~f(G,B_T)~\stackrel{\Pb}{\sim}~\nu &\mbox{and}& \Pb\left[|B_T-M|>\rho\right]~\le~\rho.
\e*
Recall that $\F= (\Fc_t)_{t \ge 0}$ is the canonical filtration generated by $B$. Set $\G:=(\Gc_t)_{t\ge 0}$ with $\Gc_t:=\sigma(\Fc_t, G)$, then $\G$ is again a Brownian filtration and more importantly, $B$ is a $\G-$Brownian motion. Take the continuous martingale $M=(M_t)_{t\ge 0}$ given by
\b*
M_t&:=&\E^{\Pb}\left[M|\Gc_t\right], \mbox{ for all } t\ge 0,
\e*
then it follows by Doob's martingale inequality that, for all $r>0$, 
\b*
&&\Pb\left[\sup_{t\le T}|B_t-M_t|>\rho^{r}\right]~\le~\frac{1}{\rho^{r}}\E^{\Pb}\left[|B_T-M|\right] \\
&\le&\frac{1}{\rho^{r}}\left\{\E^{\Pb}\left[|B_T-M|\mathbf{1}_{|B_T-M|\le \rho\}}\right]+\E^{\Pb}\left[|B_T-M|\mathbf{1}_{|B_T-M|> \rho\}}\right]\right\} \\
&\le&\frac{1}{\rho^{r}}\left\{\rho+\left(\E^{\Pb}\left[|B_T-M|^p\right]\right)^{1/p}\Pb\left[|B_T-M|> \rho\right]^{1/q}\right\} \\
&\le&\rho^{1-r}+2V^{1/p}\rho^{1/q-r}.
\e*
Notice that, in view of L\'evy's theorem, there exists a Brownian motion $W=(W_t)_{t\ge 0}$ such that
$M=\left(W_{{\langle M \rangle}_t}\right)_{t\ge 0}$ $\Pb-$almost surely. Then it follows from the time-invariance of $\Phi$ that $\Phi(M,T)=\Phi\left(W,{\langle M \rangle}_T\right)$, which yields that
\b*
&&\E^{\Pb}[\Phi(B,T)]-\pp(\nu)
~\le~\E^{\Pb}[\Phi(B,T)]-\E^{\Pb}[\Phi(M,T)] \\
&\le&\E^{\Pb}\left[\left|\Phi(B,T)-\Phi(M,T)\right|\mathbf{1}_{\{\sup_{t\le T}|B_t-M_t|\le \rho^{r}\}}\right] + 2\|\Phi\|\Pb\left[\sup_{t\le T}|B_t-M_t|> \rho^{r}\right] \\
&\le&L \rho^{r}+2\|\Phi\|\rho^{1-r}+4\|\Phi\|V^{1/p}\rho^{1/q-r}.
\e*
Optimizing with respect to $r>0$ and $\Pb\in\Pcb(\mu)$, one has a constant $C$ depending on $L$, $\|\Phi\|$ and $V$ such that 
\b*
\pp(\mu)-\pp(\nu)&\le&C\rho^{1/2q}.
\e*
Repeating the argument above by interchanging $\mu$ and $\nu$, we deduce the required result. \qed

\no Combing Propositions \ref{pro:rate} and \ref{pro:one-rate}, we obtain immediately the estimation.

\begin{Theorem}
Let Assumptions \ref{hyp:phi0} and  \ref{hyp:phi-one} hold. Then there exists a constant $C>0$ depending on $L$, $\|\Phi\|$ and $V$ such that 
\b*
0~~\le~~\pp^V(\bk^n,\bc^n)-\pp(\mu)~~\le~~C\left((\Delta \bk^n)^{1/4q}+|\bk^n|^{-p/4q^2}\right).
\e*
\end{Theorem}


\subsection{Multi-marginal case: $m\ge 1$}

For the multi-marginal case, we will make use of a duality to be specified later. The main idea here is to restrict the vector $\bl\in\Lambda^m$ in the subset $\Lambda^m_0$ defined by
\b*
\Lambda^m_0&:=&\left\{(\lambda_1,\cdots, \lambda_m)\in\Lambda^m: \lambda_i \mbox{ is } (1+L)(1+\|\Phi\|)-\mbox{Lipschitz for all } i=1,\cdots, m\right\}.
\e*
Then we have another estimation of the difference between the upper bounds of different target distributions.

\begin{Proposition}\label{pro:mul-rate}
Let Assumption \ref{hyp:phi0} hold and $p>2$. Then there exists a constant $C>0$ such that, for any centered peacock $\bnu=(\nu_1,\cdots, \nu_m)$ satisfying $\nu_m(|x|^p)\le V$, one has 
\b*
\left|\pp(\bmu)-\pp(\bnu)\right|&\le &C\sum_{i=1}^mW(\mu_i,\nu_i)^{\frac{p-2}{p-1}}.
\e*
\end{Proposition}

\noindent The strategy for the proof is to translate the embedding problem for $\bmu$ into a (modified) transport problem between the Wiener measure $\P_0$ and the target $\bmu$. To this end we equip the space $\Omb\times\R^m$ with the reward function
\b*
\xi(\om,\theta,x)&:=&
		\begin{cases}
			 \Phi(\om,\theta), ~~&\mbox{ if } \om_{\theta}=x, \\
			-\infty, ~~& \mbox{ otherwise}.
		\end{cases}
\e*
We introduce an enlarged canonical space $\Xi:=\Omb\times\R^m$ as well as its canonical elements $(B,T,X)$ defined by $B(\om,\theta,x)=\om$, $T(\om,\theta,x)=\theta$ and $X(\om,\theta,x)=x$. It is clear that $\Xi$ is a Polish space and denote by $\Qc(\Xi)$ the set of all probability measures on $\Xi$. Define further by $\Qc(\bmu)\subset\Qc(\Xi)$ the subset consisting of measures $\Q$ such that
\b*
\Q\circ (B,T)^{-1}\in\Pc,~ X_i~\stackrel{\Q}{\sim}~\mu_i \mbox{ and } \E^{\Q}[T_i]~\le~ A_i~:=~\mu_i(|x|^2) \mbox{ for all } i=1, \cdots, m.
\e*
Then the optimal transport problem is defined by
\be\label{def:pp}
\pp_{e}(\bmu)&:=&\sup_{\Q\in\Qc(\bmu)}\E^{\Q}[\xi(B,T,X)].
\ee
The problem \eqref{def:pp} is an extremal marginal problem is considered with additional constraints of the moment-type, see \textit{e.g.} Chapter 4.6.3 of Rachev \& R\"uschendorf \cite{RR} for more details. And the corresponding dual problem is defined as follows:
\be\label{def:dp}
\dd_{e}(\bmu)&:=&\inf_{(S,\bl,\bkk)\in\Dc_{e}}\bmu(\bl)+ \bkk \cdot \bmu(|x|^2),
\ee
where $\Dc_e$ denotes the set of elements $\left(S,\bl,\bkk:=(k_1,\cdots, k_m)\right)\in\Sc\times\Lambda_0^m\times\R^m$ such that 
\b*
0~\le~ k_i~\le~L, \mbox{ for all } i=1,\cdots, m
\e*
and
\b*
S_{\theta_m}+\sum_{i=1}^m\lambda_i(x_i)-k_i(\theta_i-A_i)~\ge~ \xi(\om,\theta,x), \mbox{ for all } (\om,\theta,x)\in\Omb\times\R^m.
\e*
Then Theorem 4.6.12 of \cite{RR} yields immediately the following duality.

\begin{Lemma}\label{lem:ed}
Let Assumption \ref{hyp:phi0} hold and $p\ge 2$. Then the following duality holds: 
\b*
\pp_e(\bmu)&=&\dd_e(\bmu).
\e*
\end{Lemma}

\no In view of the above duality, it remains to estimate $\mu_i(|x|^2)-\nu_i(|x|^2)$ for all $i=1,\cdots, m$, which is achieved by the next lemma.

\begin{Lemma}\label{lem:var-est}
Assume $p>2$, then there exists a constant $C>0$ such that, for any centered peacock $\bnu=(\nu_1,\cdots, \nu_m)$ satisfying $\nu_m(|x|^p)\le V$, one has 
\b*
\left|\mu_i(|x|^2)-\nu_i(|x|^2)\right|~\le ~CW(\mu,\nu)^{\frac{p-2}{p-1}}, \mbox{ for all } i=1,\cdots, m.
\e*
\end{Lemma}

\proof For the sake of simplicity we write $\mu\equiv \mu_i$ (resp. $\nu\equiv \nu_i$) for short. For all $R>0$ one has 
\b*
&&\left|\mu(|x|^2)-\nu(|x|^2)\right| \\
&\le&\left|\mu(|x|^2\mathbf{1}_{\{|x|< R\}})-\nu(|x|^2\mathbf{1}_{\{|x|< R\}})\right|+\left|\mu(|x|^2\mathbf{1}_{\{|x|\ge R\}})-\nu(|x|^2\mathbf{1}_{\{|x|\ge R\}})\right| \\
&\le&2RW(\mu,\nu)+\frac{2V}{R^{p-2}}.
\e*
Taking in particular $R=W(\mu,\nu)^{-1/(p-1)}$, one obtains some suitable constant $C$ depending only on $V$ such that
\b*
\left|\mu(|x|^2)-\nu(|x|^2)\right|&\le& CW(\mu,\nu)^{\frac{p-2}{p-1}}.
\e* 
\qed

\no Now let us turn to the proof of Proposition \ref{pro:mul-rate} by a use of Lemmas \ref{lem:ed} and \ref{lem:var-est}.

\vspace{2mm}

\noindent \textit{Proof of Proposition \ref{pro:mul-rate}.} First let us show $\pp(\bmu)= \pp_e(\bmu)$. On the one hand, one has by definition $\pp(\bmu)\le\pp_e(\bmu)$. On the other hand, for each $\Q\in\Qc(\bmu)$, it is easy to see $\E^{\Q}[\xi(B,T,X)]\neq -\infty$ if and only if $B_T\stackrel{\Q}{\sim}\bmu$, or equivalently, $\Pb:=\Q\circ (B,T)^{-1}\in\Pcb(\bmu)$. Thus
\b*
\E^{\Q}[\xi(B,T,X)]~~=~~\E^{\Pb}[\Phi(B,T)]~~\le~~\pp(\bmu),
\e*
which yields the required identity. Next, let us turn to estimate $|\pp(\bmu)-\pp(\bnu)|$. It follows by Lemma \ref{lem:ed} that there exists a sequence $(S^n, \bl^n, \bkk^n)$ such that 
\b*
\lim_{n\to\infty}\bnu(\bl^n)+\bkk^n\cdot \bnu(|x|^2)&=&\pp(\bnu),
\e*
which implies that 
\b*
\pp(\bmu)-\pp(\bnu)&=&\pp(\bmu)-\lim_{n\to\infty}\left(\bnu(\bl^n)+\bkk^n\cdot \bnu(|x|^2)\right) \\
&\le&\sup_{n\ge 1}\left\{\bmu(\bl^n)+\bkk^n\cdot \bmu(|x|^2)-\bnu(\bl^n)+\bkk^n\cdot \bnu(|x|^2)\right\} \\
&\le&\sup_{n\ge 1}\left(\bmu(\bl^n)-\bnu(\bl^n)\right)+\sup_{n\ge 1}\left(\bkk^n\cdot \bmu(|x|^2)-\bkk^n\cdot \bnu(|x|^2)\right) \\
&\le&(1+L)(1+\|\Phi\|)\sum_{i=1}^mW(\mu_i,\nu_i)+L\sum_{i=1}^m\left|\mu_i(|x|^2)-\nu_i(|x|^2)\right| \\
&\le& C\sum_{i=1}^mW(\mu_i,\nu_i)^{\frac{p-2}{p-1}},
\e*
where $C$ depends only on $L$, $\|\Phi\|$ and $V$. Repeating the reasoning above by interchanging $\bmu$ and $\bnu$ we get the required result. \qed 

\vspace{2mm}

\no Combing Propositions \ref{pro:rate} and \ref{pro:mul-rate}, we obtain immediately the main result under the following assumption.

\begin{Assumption} \label{hyp:stab2}
The sequence $(\bk^n)_{n\ge 1}$ satisfies 
\b*
\lim_{n\to\infty}|\bk^n|~=~+\infty&\mbox{and}&\lim_{n\to\infty} |\bk^n |\sqrt{\Delta \bk^n}~=~0.
 \e*
\end{Assumption} 

\begin{Theorem}
Let Assumption \ref{hyp:phi0} hold and $p>3$. Then there exists a constant $C>0$ depending on $L$, $\|\Phi\|$ and $V$ such that 
\b*
0~~\le~~\pp^V(\bk^n,\bc^n)-\pp(\bmu)~~\le~~C|\bk^n|^{\frac{p-2}{p-1}}\left(\sqrt{\Delta \bk^n}+|\bk^n|^{-p/2q}\right)^{\frac{p-2}{p-1}}.
\e*
\end{Theorem}

\end{document}